\documentclass[12pt]{article}
\usepackage[]{amsmath,amssymb}
\usepackage{amscd}
\usepackage{latexsym}
\usepackage{cite}

\newtheorem{definition}{Definition}[section]
\newtheorem{theorem}[definition]{Theorem}
\newtheorem{lemma}[definition]{Lemma}
\newtheorem{corollary}[definition]{Corollary}

\newtheorem{example}[definition]{Example}

\newtheorem{proposition}[definition]{Proposition}
\newtheorem{notation}[definition]{Notation}

\def\F{\mathbb F}

\begin{document}
\title{\bf An algebra generated by two sets of\\
mutually orthogonal idempotents}
\author{
Tatsuro Ito\footnote{Supported in part by JSPS grant
18340022.} $\;$   and
Paul Terwilliger}
\date{}

\maketitle
\begin{abstract}
\noindent For a field $\F$ and
an integer $d\geq 1$, we consider the universal
associative 
$\F$-algebra $A$ generated by
two sets of $d+1$ mutually orthogonal idempotents.
We display four bases for the $\F$-vector
space $A$ that we find attractive. We determine
how these bases are related to each other.
We describe how the multiplication
in $A$ looks with respect to our bases.
Using our bases we obtain
an infinite nested sequence of
2-sided ideals for $A$.
Using our bases we obtain an infinite
exact sequence involving a certain 
$\F$-linear map
$\partial: A \to A$.
We obtain several results concerning
the kernel of $\partial$; 
for instance this kernel
 is 
a subalgebra of $A$ that is free of
rank $d$.




\bigskip
\noindent
{\bf Keywords}. 
Coproduct, free algebra, nonrepeating word. 
 \hfil\break
\noindent {\bf 2000 Mathematics Subject Classification}. 
Primary: 16S10. 
 \end{abstract}

\section{The algebra $A$}

\noindent Throughout the paper $\F$ denotes a field.
All unadorned tensor products are meant to be over $\F$.
An algebra is meant to be associative and have a 1.

\medskip
\noindent
We now introduce our topic.

\begin{definition}
\label{def:t}
\rm
Let $d$ denote a positive integer.
Let $A=A(d,\F)$ denote the $\F$-algebra 
defined by
generators $\lbrace e_i\rbrace_{i=0}^d$,
 $\lbrace e^*_i\rbrace_{i=0}^d$ and the
following relations:
\begin{eqnarray}
&&e_ie_j = \delta_{i,j} e_i,
\qquad \qquad
e^*_ie^*_j = \delta_{i,j} e^*_i,
\qquad \qquad (0 \leq i,j\leq d),
\label{eq:rel1}
\\
&&1 = \sum_{i=0}^d e_i, 
\qquad \qquad 
1 = \sum_{i=0}^d e^*_i.
\label{eq:rel2}
\end{eqnarray}
Here $\delta_{i,j}$ denotes the Kronecker delta.
\end{definition}

\begin{definition}
\label{def:gens}
\rm
Referring to Definition \ref{def:t},
we call 
 $\lbrace e_i\rbrace_{i=0}^d$ and
 $\lbrace e^*_i\rbrace_{i=0}^d$ the {\it idempotent 
generators} for $A$.
We say that the  
 $\lbrace e^*_i\rbrace_{i=0}^d$ are {\it starred}
and the 
 $\lbrace e_i\rbrace_{i=0}^d$ are {\it nonstarred}.
\end{definition}

\noindent We now briefly explain how $A$ can be viewed as a
coproduct in the sense of 
Bergman\cite{bergman1,bergman2}.
As we will see in Theorem
\ref{thm:tbasis}, the elements
$\lbrace e_i\rbrace_{i=0}^d$ are linearly
independent in $A$ and hence form a basis
for a subalgebra of $A$ 
denoted
 $A_1$. Similarly
the elements
$\lbrace e^*_i\rbrace_{i=0}^d$ 
form a basis for a subalgebra of 
$A$ denoted $A^*_1$.
By construction $A$ is generated by
$A_1, A^*_1$. The $\F$-algebras 
$A_1$ and $A^*_1$
are each isomorphic to a direct sum
of $d+1$ copies of $\F$.
The elements
$\lbrace e_i\rbrace_{i=0}^d$
(resp. 
$\lbrace e^*_i\rbrace_{i=0}^d$)
are the primitive idempotents of
$A_1$ (resp. $A^*_1$).
Since no relation in 
(\ref{eq:rel1}), 
(\ref{eq:rel2}) 
involves both $A_1$ and $A^*_1$,
the algebra $A$ is the coproduct of
$A_1$ and $A^*_1$ in the sense of
Bergman\cite[Section~1]{bergman1}.
As part of his comprehensive study of
coproducts, Bergman determined the
units and zero-divisors in $A$
\cite[Corollary 2.16]{bergman1}.

\medskip
\noindent 
Our goal in this article is to describe four
bases
for $A$ that we find attractive.
We  determine how these bases are related
to each other. We describe how the multiplication
in $A$ looks with respect to these bases. 
Using our bases we obtain an infinite
nested sequence of
2-sided ideals for $A$.
Using our bases we obtain an infinite
exact sequence involving a certain 
$\F$-linear map $\partial: A \to A$.
We show that the kernel $F$ of $\partial$
is a subalgebra of $A$ that is free
of rank $d$. 
We show that $F$ is generated by the elements
$\lbrace e_i-e^*_i\rbrace_{i=1}^d$.
We show that each of the $\F$-linear maps
\begin{eqnarray*}
&&F \otimes A_1 \;\to\; A 
\qquad \qquad \qquad 
F \otimes A^*_1 \; \;\to \;\; A
\\
&&u \otimes v \; \; \mapsto \; \; uv
\qquad \qquad \qquad 
u \otimes v \;  \;\mapsto \;  \; uv
\end{eqnarray*}
is an isomorphism of $\F$-vector spaces.
We will define our bases after a few comments.

\medskip
\noindent The following three lemmas are about symmetries of $A$;
their proofs are routine and left to the reader.

\begin{lemma}
\label{lem:aut}
There exists a unique $\F$-algebra automorphism
of $A$ that sends
\begin{eqnarray*}
e_i \mapsto e^*_i, 
\qquad \qquad e^*_i \mapsto e_i
\end{eqnarray*}
for $0 \leq i \leq d$.
Denoting this automorphism by $*$ 
we have $x^{**}=x$ for all $x \in A$.
\end{lemma}

\noindent By an  
{\it $\F$-algebra antiautomorphism} of $A$
we mean an isomorphism of $\F$-vector spaces
$\rho: A\to A$ such that $(xy)^\rho=y^\rho x^\rho$
for all $x,y \in A$.
\begin{lemma}
\label{lem:antiaut}
There exists a unique $\F$-algebra antiautomorphism
$\dagger $  of $A$ that fixes each idempotent generator.
We have
$x^{\dagger \dagger}=x$ for all $x \in A$.
\end{lemma}

\begin{lemma}
\label{lem:mapscom}
The maps $*$ and $\dagger$ commute.
\end{lemma}

\noindent Let $X$ denote a subset of $A$.
By the {\it relatives} of $X$ we mean
the subsets
$X$, $X^*$, $X^\dagger$,
$X^{*\dagger}$.

\section{Four bases for the vector space $A$}

In this section we display four bases for the
$\F$-vector space $A$.

\begin{definition}
\label{def:word}
\rm
A pair of idempotent generators for $A$ is called
{\it alternating}
whenever one of them is starred and the other is
nonstarred.
For an integer $n\geq 1$, by a {\it word of length $n$} in $A$
we mean a product $g_1g_2 \cdots g_n$ such that
 $\lbrace g_i\rbrace_{i=1}^n$
are idempotent generators for $A$ and $ g_{i-1}, g_i$
are alternating for $2 \leq i \leq n$.
The  word 
$g_1g_2 \cdots g_n$ 
is said to {\it begin} with $g_1$
and {\it end} with $g_n$. 
\end{definition}

\begin{example}
\rm For $d=2$ we display the words in $A$ that have length
3 and begin with $e_0$.
\begin{eqnarray*}
&&e_0e^*_0e_0, \qquad e_0e^*_0e_1, \qquad e_0e^*_0e_2,
\\
&&e_0e^*_1e_0, \qquad e_0e^*_1e_1, \qquad e_0e^*_1e_2,
\\
&&e_0e^*_2e_0, \qquad e_0e^*_2e_1, \qquad e_0e^*_2e_2.
\end{eqnarray*}
\end{example}

\begin{definition}
\rm
For an idempotent generator $e_i$ or $e^*_i$ we call
$i$ the {\it index} of the generator. 
A word $g_1 g_2 \cdots g_n$ in $A$   
is called {\it nonrepeating} (or {\it NR})
whenever $g_{j-1}, g_j$ do not have the same index for
 $2 \leq j \leq n$.
\end{definition}

\begin{example} \rm
For $d=2$ we display the NR words in $A$ that have length 3
and begin with $e_0$.
\begin{eqnarray*}
e_0e^*_1e_0, \qquad
e_0e^*_1e_2, \qquad
e_0e^*_2e_0, \qquad
e_0e^*_2e_1.
\end{eqnarray*}
\end{example}


\begin{theorem}
\label{thm:tbasis}
Each of the following is a basis for the $\F$-vector space
$A$:
\begin{itemize}
\item[\rm (i)]
The set of  NR words in $A$ that end with a nonstarred element.
\item[\rm (ii)]
The set of  NR words in $A$ that end with a starred element.
\item[\rm (iii)]
The set of NR words in $A$ that begin with a nonstarred element.
\item[\rm (iv)]
The set of NR words in $A$ that begin with a starred element.
\end{itemize}
\end{theorem}
\noindent {\it Proof:} 
(i)
Let $S$ denote the set of
  NR words in $A$ that end with a nonstarred element.
We first show that $S$ spans $A$.
Let $A'$ denote the subspace of $A$ spanned
by $S$.
To obtain $A'=A$ 
it suffices to show that $A'$ is a left ideal of $A$
that contains $1$.
To show that $A'$ is a left ideal of $A$,
it suffices to show that  
$e_iA' \subseteq A'$ and $e^*_iA' \subseteq A'$
for $0 \leq i \leq d$. 
For a word $w=g_1g_2\cdots g_n$
 in $S$ and  $0 \leq i \leq d$ 
 we show that each of $e_iw$, $e^*_iw$ is contained in $A'$. 
Let $j$ denote the index of $g_1$.
Invoking
(\ref{eq:rel2})
we may assume without loss that $i \not=j$.
First assume $n$ is odd, so that  
$g_1=e_j$.
Since $e_ie_j=0$ we have
$e_iw=0$, so $ e_iw\in A'$.
Also $e^*_iw=e^*_ig_1g_2\cdots g_n$  
is a word in $S$, so $e^*_iw \in A'$.
Next assume $n$ is even, so that  $g_1=e^*_j$.
Then $e_iw=e_ig_1g_2\cdots g_n$  
is a word in $S$, so $e_iw \in A'$.
Since $e^*_ie^*_j=0$ we have
$e^*_iw=0$, so $e^*_iw\in A'$.
We have shown 
$A'$ is a left 
ideal of $A$.
The ideal $A'$ contains 1, since
 $e_i \in S$ for $0 \leq i \leq d$ and
$1=\sum_{i=0}^d e_i$.
We have shown $A'$ is a left ideal of $A$ that contains $1$, so $A'=A$.
Therefore $S$ spans $A$.
Next we show that
the elements of $S$ are linearly independent.
Let $\cal S$ denote the set of sequences
$(r_1, r_2, \ldots, r_n)$
such that (i) $n$ is a positive integer; (ii)
each of $r_1, r_2, \ldots, r_n$ is contained in  
the set $\lbrace 0,1,\ldots, d\rbrace$;
(iii) $r_{i-1}\not=r_i$ for $2 \leq i \leq n$.
Let $V$ denote the vector space over $\F$ consisting
of those formal linear combinations of $ \cal S$
that have finitely many nonzero coefficients.
The set $\cal S$ is a basis for $V$.
For $0 \leq i \leq d$ we  define linear transformations
$E_i: V \to V$ and 
$E^*_i: V \to V$. To this end
we give the actions of $E_i$ and $E^*_i$ on 
 $\cal S$. Pick an element $(r_1, r_2, \ldots, r_n) \in \cal S$.  
The actions of $E_i$ and $E^*_i$ on 
  $(r_1, r_2, \ldots, r_n)$
are given in
the table below.

\begin{center}
\begin{tabular}{c|c|c}
Case  &  $E_i.(r_1, \ldots, r_n)$ & $E^*_i.(r_1, \ldots, r_n)$ \\
\hline
$r_1=i$,  $\;n$ odd &
$(r_1, \ldots, r_n)$ & $(r_1, \ldots, r_n)
 - 
\sum_{\stackrel{0 \leq j \leq d}{j \neq i}}
 (j, r_1,  \ldots, r_n)$
\\
$r_1\not=i$, $\;n$ odd & $0$ & $(i,r_1,  \ldots, r_n)$ 
\\
$r_1=i$, $\;n$ even  & $(r_1,  \ldots, r_n)-
\sum_{\stackrel{0 \leq j \leq d}{j \neq i}}
 (j, r_1,  \ldots, r_n)$ & $(r_1, \ldots, r_n)$
 \\
$r_1\not=i$, $\;n$ even & $(i,r_1,  \ldots, r_n)$ & $0$ 
\end{tabular}
\end{center}

\bigskip
\noindent Using the table,
\begin{eqnarray}
&&E_iE_j = \delta_{i,j} E_i,
\qquad \qquad
E^*_iE^*_j = \delta_{i,j} E^*_i,
\qquad \qquad (0 \leq i,j\leq d),
\label{eq:Erel1}
\\
&&1 = \sum_{i=0}^d E_i, 
\qquad \qquad 
1 = \sum_{i=0}^d E^*_i.
\label{eq:Erel2}
\end{eqnarray}
Comparing
(\ref{eq:Erel1}),
(\ref{eq:Erel2})
with
(\ref{eq:rel1}), 
(\ref{eq:rel2}) we find that $V$ has an $A$-module
structure such that $e_i$ (resp. $e^*_i$)
acts on $V$ as $E_i$ (resp. $E^*_i$) for $0 \leq i \leq d$.
Define the element 
$\Delta \in V$ by
$\Delta= \sum_{i=0}^d (i)$, and consider the
linear transformation
$\gamma: A \to V$ that sends
$x \mapsto x.\Delta$ for all $x \in A$.
For each word 
$w=g_1g_2\cdots g_n$ in $S$
we find 
$\gamma(w)= ({\overline {g_1}} ,
 {\overline {g_2}}, \ldots,   
{\overline {g_n}} )$ where 
$\overline g$ denotes the index of $g$.
Thus the restriction of  $\gamma$ to $S$ gives
a bijection $S \to {\cal S}$. 
The elements of $\cal S$ are linearly 
independent and $\gamma$ is linear,
so the elements of $S$ are
linearly
independent.
We have shown $S$ is a basis for $A$.
\\
\noindent (ii) Apply the automorphism
$*$ to the basis in (i) above.
\\
\noindent (iii), (iv) Apply the antiautomorphim
$\dagger$ to the bases in (i), (ii) above.
\hfill $\Box$ \\

\section{How the four bases for $A$ are related}

\noindent In this section we obtain some identities
that effectively give the transition matrix
between any two  bases from Theorem
\ref{thm:tbasis}. 

\medskip

\begin{notation}
\label{not:index}
\rm
Let $w=g_1g_2\cdots g_n$ denote a 
  word in $A$, with
$g_n$ nonstarred.
We represent $w$
by the sequence $(r_1, r_2, \ldots, r_n)$, where
$r_j$ denotes the index of $g_j$ for $1 \leq j \leq n$. 
We represent $w^*$ by
 $(r_1, r_2, \ldots, r_n)^*$.
\end{notation}

\begin{example}
\rm
We display some  words in $A$ along
with their
notation.

\begin{center}
\begin{tabular}{c|c}
word  &  notation \\
\hline
$e_0e^*_2e_1$ & $(0,2,1)$
\\
$e^*_1e_0e^*_2e_1$ & $(1,0,2,1)$
\\
$e^*_0e_2e^*_1$ & $(0,2,1)^*$
\\
$e_1e^*_0e_2e^*_1$ & $(1,0,2,1)^*$ 
\end{tabular}
\end{center}
\end{example}

\noindent The next result effectively gives
the transition matrix between any two
bases from Theorem \ref{thm:tbasis}.

\begin{theorem}
\label{thm:trans}
With reference to Notation 
\ref{not:index}, and
for each basis vector $(r_1, r_2, \ldots, r_n)$ from
Theorem
\ref{thm:tbasis}(i), the element 
\begin{eqnarray*}
(r_1, r_2, \ldots, r_n)
\quad + \quad 
(-1)^n(r_1, r_2, \ldots, r_n)^*
\end{eqnarray*}
is equal to
\begin{eqnarray*}
&& 
 \sum_{\stackrel{0 \leq j \leq d}{j \neq r_1}}
(j,r_1, r_2, \ldots, r_n)
\quad +  \quad 
 \sum_{\ell=1}^{n-1} (-1)^{\ell}
\sum_{\stackrel{0 \leq j \leq d}{j \neq r_{\ell},\; j\neq r_{\ell+1}}}
(r_1, r_2, \ldots, r_{\ell}, j, r_{\ell+1}, \ldots, r_n)
\\
&&  \qquad \quad 
+\quad  
(-1)^n\sum_{\stackrel{0 \leq j \leq d}{j \neq r_n}}
( r_1, r_2, \ldots, r_n,j),
\end{eqnarray*}
\noindent and also equal to
\begin{eqnarray*}
&& 
 \sum_{\stackrel{0 \leq j \leq d}{j \neq r_n}}
(r_1, r_2, \ldots, r_n,j)^*
\quad +  \quad 
 \sum_{\ell=1}^{n-1} (-1)^{n-\ell}
\sum_{\stackrel{0 \leq j \leq d}{j \neq r_{\ell}, \;j\neq r_{\ell+1}}}
(r_1, r_2, \ldots, r_{\ell}, j, r_{\ell+1}, \ldots, r_n)^*
\\
&&  \qquad \quad 
+\quad  (-1)^n 
\sum_{\stackrel{0 \leq j \leq d}{j \neq r_1}}
(j, r_1, r_2, \ldots, r_n)^*.
\end{eqnarray*}
\end{theorem}
\noindent {\it Proof:} 
To obtain the first assertion, define
\begin{eqnarray}
\phi_0  &=& 
 \sum_{\stackrel{0 \leq j \leq d}{j \neq r_1}}
(j,r_1, r_2, \ldots, r_n),
\label{eq:phi0}
\\
 \phi_{\ell} &=& 
\sum_{\stackrel{0 \leq j \leq d}{j \neq r_{\ell},\; j\neq r_{\ell+1}}}
(r_1, r_2, \ldots, r_{\ell}, j, r_{\ell+1}, \ldots, r_n)
\qquad \qquad 
(1 \leq \ell \leq n-1),
\label{eq:phil}
\\
\phi_n
&=& 
\sum_{\stackrel{0 \leq j \leq d}{j \neq r_n}}
( r_1, r_2, \ldots, r_n,j).
\label{eq:phin}
\end{eqnarray}
Evaluating 
(\ref{eq:phi0})--(\ref{eq:phin})
using 
(\ref{eq:rel2}) we find
\begin{eqnarray}
\phi_0 &=& (r_1, r_2, \ldots, r_n)-
 (r_1,r_1, r_2, \ldots, r_n),
\label{eq:phi0ev}
\\
\phi_{\ell} &=& -(r_1, r_2, \ldots, r_{\ell}, r_{\ell},r_{\ell+1},\ldots, r_n)-
(r_1, r_2, \ldots, r_{\ell}, r_{\ell+1},r_{\ell+1},\ldots, r_n),
\label{eq:philev}
\\
\phi_n &=& (r_1, r_2, \ldots, r_n)^*-
 (r_1, r_2, \ldots, r_n, r_n).
\label{eq:phinev}
\end{eqnarray}
Combining
(\ref{eq:phi0ev})--(\ref{eq:phinev}) we obtain
\begin{eqnarray*}
\phi_0 + \sum_{\ell=1}^{n-1} (-1)^{\ell} \phi_{\ell}
+ (-1)^n \phi_n &=& 
(r_1, r_2, \ldots, r_n)+
(-1)^n(r_1, r_2, \ldots, r_n)^*,
\end{eqnarray*}
and the first assertion follows. The second assertion is similarly
obtained.
\hfill $\Box$ \\

\begin{example}
\rm
Assume $d=2$. Then for $n=1$ and $r_1=1$ the assertions
of Theorem
\ref{thm:trans} become
\begin{eqnarray*}
e_1 - e^*_1 
&=& e^*_0e_1+e^*_2e_1
 -e^*_1e_0 - e^*_1e_2 
\\
&=& e_1e^*_0 + e_1e^*_2-e_0e^*_1-e_2e^*_1.
\end{eqnarray*}
For $n=2$ and $(r_1, r_2)=(1,0)$ 
 the assertions
of Theorem
\ref{thm:trans} become
\begin{eqnarray*}
e^*_1e_0+e_1e^*_0 &=& 
 e_0e^*_1e_0+
e_2e^*_1e_0
-e_1 e^*_2  e_0
+e_1e^*_0e_1+
e_1e^*_0e_2
\\
&=&
e^*_1e_0e^*_1+
e^*_1e_0e^*_2
-e^*_1 e_2 e^*_0
+ e^*_0e_1e^*_0
+ e^*_2e_1e^*_0.
\end{eqnarray*}
\end{example}

\section{The product of basis elements}

Consider the basis for $A$ from
Theorem
\ref{thm:tbasis}(i).
We now take two elements from this basis,
and write the product
as a linear combination of elements from the basis.

\begin{theorem}
\label{thm:prod}
Let $(r_1, r_2, \ldots, r_n)$
and 
 $(r'_1, r'_2, \ldots, r'_m)$
denote basis vectors from Theorem
\ref{thm:tbasis}(i). 
Then the product
\begin{eqnarray}
(r_1, r_2, \ldots, r_n)\cdot 
 (r'_1, r'_2, \ldots, r'_m)
\label{eq:prod}
\end{eqnarray}
is the following linear combination of basis
vectors from Theorem
\ref{thm:tbasis}(i). 
\begin{itemize}
\item[\rm (i)]
Assume $m$ is odd and $r_n \neq r'_1$. Then
(\ref{eq:prod}) is zero.
\item[\rm (ii)]
Assume $m$ is odd and $r_n=r'_1$. Then
(\ref{eq:prod}) is equal to
\begin{eqnarray*} 
(r_1, r_2, \ldots, r_n,
 r'_2, \ldots, r'_m).
\end{eqnarray*}
\item[\rm (iii)]
Assume $m$ is even and $r_n \neq r'_1$. Then
(\ref{eq:prod}) is equal to 
\begin{eqnarray*}
(r_1, r_2, \ldots, r_n,
 r'_1,r'_2, \ldots, r'_m).
\end{eqnarray*}
\item[\rm (iv)]
Assume $m$ is even and $r_n = r'_1$. Then
(\ref{eq:prod}) is equal to 
\begin{eqnarray*}
&&(-1)^{n+1}(r_1, r_2, \ldots, r_n,r'_2, \ldots, r'_m)\quad + \quad  (-1)^n \sum_{\stackrel{0 \leq j \leq d}{j \neq r_1}}
 (j, r_1, r_2, \ldots, r_n, r'_2, \ldots, r'_m)
\\
&&\quad +\quad \sum_{\ell=1}^{n-1} (-1)^{n-\ell}
 \sum_{\stackrel{0 \leq j \leq d}{j \neq r_{\ell}, \;j\neq r_{\ell+1}}}
 (r_1,r_2, \ldots, r_{\ell},j,r_{\ell+1},\ldots, r_n, r'_2, \ldots, r'_m).
\end{eqnarray*}
\end{itemize}
\end{theorem}
\noindent {\it Proof:} 
(i)--(iii) Routine.
\\
\noindent (iv) In line 
(\ref{eq:prod}), evaluate 
$(r_1, r_2, \ldots, r_n)$
using the second identity in
Theorem
\ref{thm:trans}, and simplify the result.
\hfill $\Box$ \\

\noindent 
Now consider 
the basis for $A$
from Theorem
\ref{thm:tbasis}(ii),
and 
the basis for $A$ from
Theorem
\ref{thm:tbasis}(i).
In the next result, we take an element from
the first basis and an element
from the second basis, and
write the product as a 
 linear
combination of elements from the
second basis.

\begin{theorem}
\label{thm:prod2}
In the notation of Theorem
\ref{thm:prod},
 the product
\begin{eqnarray}
(r_1, r_2, \ldots, r_n)^*\cdot 
 (r'_1, r'_2, \ldots, r'_m)
\label{eq:prod2}
\end{eqnarray}
is the following linear combination of basis
vectors from Theorem
\ref{thm:tbasis}(i). 
\begin{itemize}
\item[\rm (i)]
Assume $m$ is even and $r_n \neq r'_1$. Then
(\ref{eq:prod2}) is $0$.
\item[\rm (ii)]
Assume $m$ is even and $r_n=r'_1$. Then
(\ref{eq:prod2}) is equal to
\begin{eqnarray*} 
(r_1, r_2, \ldots, r_n,
 r'_2, \ldots, r'_m).
\end{eqnarray*}
\item[\rm (iii)]
Assume $m$ is odd and $r_n \neq r'_1$. Then
(\ref{eq:prod2}) is equal to 
\begin{eqnarray*}
(r_1, r_2, \ldots, r_n,
 r'_1,r'_2, \ldots, r'_m).
\end{eqnarray*}
\item[\rm (iv)]
Assume $m$ is odd and $r_n = r'_1$. Then
(\ref{eq:prod2}) is equal to 
\begin{eqnarray*}
&&(-1)^{n+1}(r_1, r_2, \ldots, r_n,r'_2, \ldots, r'_m)\quad + \quad  (-1)^n \sum_{\stackrel{0 \leq j \leq d}{j \neq r_1}}
 (j, r_1, r_2, \ldots, r_n, r'_2, \ldots, r'_m)
\\
&&\quad +\quad \sum_{\ell=1}^{n-1} (-1)^{n-\ell}
 \sum_{\stackrel{0 \leq j \leq d}{j \neq r_{\ell}, \;j\neq r_{\ell+1}}}
 (r_1,r_2, \ldots, r_{\ell},j,r_{\ell+1},\ldots, r_n, r'_2, \ldots, r'_m).
\end{eqnarray*}
\end{itemize}
\end{theorem}
\noindent {\it Proof:} 
(i)--(iii) Routine.
\\
\noindent (iv) In line 
(\ref{eq:prod2}), evaluate 
$(r_1, r_2, \ldots, r_n)^*$
using the first identity in
Theorem
\ref{thm:trans}, and simplify the result.
\hfill $\Box$ \\

\section{The subspaces $A_n$}

\noindent
 In this section we introduce
some subspaces $A_n$ of $A$, and use them to
 interpret
our results so far.  

\begin{definition}
\label{def:tn}
\rm
For an integer $n\geq 1$ let $A_n$ denote the subspace of
$A$ spanned by the NR words that have length $n$
and end with a nonstarred element.
\end{definition}

\begin{lemma}
\label{lem:sd}
For $n\geq 1$ we display a basis for each
relative of $A_n$.

\begin{center}
\begin{tabular}{c|c}
space  &  basis \\
\hline
$A_n$ &
the NR words in $A$ that have length $n$ and end with a nonstarred element
\\
$A^*_n$ &
the NR words in $A$ that have length $n$ and end with a starred element
\\
$A^\dagger_n$ &
the NR words in $A$ that have length $n$ and begin with a nonstarred element
\\
$A^{*\dagger}_n$ &
the NR words in $A$ that have length $n$ and begin with a starred element
\end{tabular}
\end{center}
\end{lemma}
\noindent {\it Proof:} 
Immediate from Lemma
\ref{lem:aut}, 
Lemma
\ref{lem:antiaut}, 
and Theorem
\ref{thm:tbasis}.
\hfill $\Box$ \\

\begin{lemma}
\label{lem:dimtn}
For $n\geq 1$ each relative of 
$A_n$
has dimension $(d+1)d^{n-1}$.
\end{lemma} 
\noindent {\it Proof:} 
Apply Lemma
\ref{lem:sd} and a routine counting argument.
\hfill $\Box$ \\

\begin{lemma}
\label{lem:eveodd}
The following {\rm (i), (ii)} hold
for all integers $n\geq 1$.
\begin{itemize}
\item[\rm (i)] Suppose $n$ is even. Then $A^{*\dagger}_n=A_n$ and
$A^*_n=A^\dagger_n$.
\item[\rm (ii)] Suppose $n$ is odd. Then $A^{*\dagger}_n=A^*_n$ and
$A^\dagger_n=A_n$.
\end{itemize}
\end{lemma}
\noindent {\it Proof:} Pick a word $w$ in $A$ of length $n$.
If $n$ is even, then $w$ 
begins with a starred element if and only if $w$ ends
with a nonstarred element. If $n$ is odd,
then $w$ begins with a starred element if and only if
$w$ ends with a starred element. The result follows.
\hfill $\Box$ \\

\begin{theorem}
\label{lem:ds}
Each of the following sums is direct.
\begin{eqnarray*}
&&A = \sum_{n=1}^\infty A_n,
\qquad \qquad 
A = \sum_{n=1}^\infty A^*_n,
\\
&&
A = \sum_{n=1}^\infty A^\dagger_n,
\qquad \qquad 
A = \sum_{n=1}^\infty A^{*\dagger}_n.
\end{eqnarray*}
\end{theorem} 
\noindent {\it Proof:} Combine
Theorem
\ref{thm:tbasis}
and Lemma \ref{lem:sd}.
\hfill $\Box$ \\

\begin{theorem}
\label{thm:interp2}
For  $n\geq 1$ and $x \in A_n$,
\begin{eqnarray*}
x +(-1)^n x^* \in A_{n+1}\cap A^*_{n+1}.
\end{eqnarray*}
\end{theorem}
\noindent {\it Proof:} 
By Definition
\ref{def:tn} we may assume without loss that
$x$
is an NR word in $A$ that has length $n$ and ends with
a nonstarred element. Now 
$x+(-1)^n x^* \in A_{n+1}$
by the first assertion
of Theorem
\ref{thm:trans}, and
$x+(-1)^n x^* \in A^*_{n+1}$
by the second assertion of
 Theorem
\ref{thm:trans}.
The result follows.
\hfill $\Box$ \\

\begin{corollary}
\label{cor:interp2}
For  $n\geq 1$,
\begin{eqnarray*}
A_n +  
 A_{n+1}\cap A^*_{n+1}
= 
A^*_n +  
A_{n+1}\cap A^*_{n+1}.
\end{eqnarray*}
\end{corollary}
\noindent {\it Proof:} 
This is a routine consequence of 
Theorem
\ref{thm:interp2}.
\hfill $\Box$ \\

\noindent For subsets $X,Y$ of $A$ 
let $XY$ denote the subspace of $A$
spanned by $\lbrace xy \,|\, x \in X, \;y \in Y\rbrace $.

\begin{theorem}
\label{thm:ptntm}
For positive integers $n,m$ the products $A_n A_m$
and $A^*_n A_m$ 
are described as follows.
\begin{itemize}
\item[\rm (i)]
 Assume $m$ is odd. Then
\begin{eqnarray}
 A_n A_m \subseteq A_{n+m-1},
\qquad \qquad 
 A^*_n A_m \subseteq A_{n+m}+A_{n+m-1}.
\label{eq:tntmodd}
\end{eqnarray}
\item[\rm (ii)]
Assume $m$ is even. Then
\begin{eqnarray}
 A_n A_m \subseteq A_{n+m}+A_{n+m-1},
\qquad \qquad 
 A^*_n A_m \subseteq A_{n+m-1}.
\label{eq:tntmeven}
\end{eqnarray}
\end{itemize}
\end{theorem}
\noindent {\it Proof:} 
In (\ref{eq:tntmodd})
and (\ref{eq:tntmeven}) the
inclusions on the left
follow from Theorem
\ref{thm:prod},
and the inclusions on the right
follow from
Theorem \ref{thm:prod2}.
\hfill $\Box$ \\

\noindent In Section 9 we will obtain a more
detailed version of Theorem 
\ref{thm:ptntm}.

\section{The ideals $A_{\geq n}$}

\noindent Motivated by
Corollary
\ref{cor:interp2}
and
 Theorem
\ref{thm:ptntm} we consider the following subspaces of $A$.

\begin{definition}
\label{def:tgeqn}
\rm
For  $n\geq 1$ define
\begin{eqnarray*}
A_{\geq n} = A_n + A_{n+1} + \cdots
\end{eqnarray*}
\end{definition}

\begin{theorem} For  $n\geq 1$
the space $A_{\geq n}$
is a 2-sided ideal of $A$.
\end{theorem}
\noindent {\it Proof:} 
This is a routine consequence of the
inclusions on the left in
(\ref{eq:tntmodd}),
(\ref{eq:tntmeven}).
\hfill $\Box$ \\

\begin{theorem}
\label{thm:same}
 For $n\geq 1$
we have
\begin{eqnarray*}
 A^*_{\geq n}=A_{\geq n},  \qquad \qquad
 A^\dagger_{\geq n}=A_{\geq n}.
\end{eqnarray*}
\end{theorem}
\noindent {\it Proof:} 
For $m\geq 1$ we obtain 
$A_m \subseteq A^*_m+A^*_{m+1}$
and 
$A^*_m \subseteq A_m+A_{m+1}$
from Corollary
\ref{cor:interp2}.
Therefore 
 $A^*_{\geq n}=A_{\geq n}$.
For $m\geq 1$ the space $A^\dagger_m$
is one of $A_m$, $A^*_m$ by
Lemma
\ref{lem:eveodd},
and each of $A_m$, $A^*_m$ is contained in
$A_m+A_{m+1}$, 
so $A^\dagger_m \subseteq A_m+A_{m+1}$.
In this inclusion we apply $\dagger$ to each side
and find
 $A_m \subseteq A^\dagger_m+A^\dagger_{m+1}$.
Therefore
 $A^\dagger_{\geq n}=A_{\geq n}$.
\hfill $\Box$ \\

\begin{lemma}
\label{lem:prodideal}
For positive integers $n,m$
the product $A_{\geq n}
A_{\geq m}$ is contained in
$A_{\geq n+m-1}$.
\end{lemma}
\noindent {\it Proof:} 
This follows from
Definition \ref{def:tgeqn}
and the products on the left
in (\ref{eq:tntmodd}),
 (\ref{eq:tntmeven}).
\hfill $\Box$ \\

\section{The map $\partial:A\to A$}

\noindent Motivated by Theorem
\ref{thm:interp2} we consider the following map.

\begin{lemma}
\label{lem:map}
\rm
There exists a unique $\F$-linear transformation 
 $\partial:A\to A$ 
 such that for $n\geq 1$,
\begin{eqnarray}
\partial(x)= x+(-1)^n x^*   
\qquad \qquad                 (\forall x \in A_n).
\label{eq:par}
\end{eqnarray}
\end{lemma}
\noindent {\it Proof:} 
By Theorem
\ref{lem:ds} the sum
$A=\sum_{n=1}^\infty A_n$
is direct.
\hfill $\Box$ \\

\begin{lemma}
\label{lem:parset}
With reference to Lemma
\ref{lem:map}
we have $\partial (A_n)\subseteq A_{n+1}$  for $n \geq 1$.
\end{lemma}
\noindent {\it Proof:} 
Immediate from Theorem
\ref{thm:interp2}.
\hfill $\Box$ \\

\begin{lemma}
\label{lem:bound}
With reference to Lemma
\ref{lem:map} the following
{\rm (i), (ii)} hold for all $x \in A$.
\begin{itemize}
\item[\rm (i)]
$\partial (\partial(x))=0$.
\item[\rm (ii)] 
$\partial(x^*)=-(\partial(x))^*$.
\end{itemize}
\end{lemma}
\noindent {\it Proof:} Without loss we may assume
$x \in A_n$ for some
$n\geq 1$. \\
\noindent (i)
Observe $\partial(x)\in A_{n+1}$ by 
Lemma
\ref{lem:parset}, so 
$\partial(\partial(x))=\partial(x)-(-1)^n(\partial(x))^*$
by
(\ref{eq:par}).
In line
(\ref{eq:par}) 
we apply $*$ to both sides and get
$(\partial(x))^*=(-1)^n \partial(x)$.
The result follows.
\\
\noindent (ii)
In line
(\ref{eq:par}) 
we apply $\partial$ to both sides and use
(i) above to get
$\partial(x^*)=(-1)^{n-1} \partial(x)$.
We observed
$(\partial(x))^*=(-1)^n \partial(x)$ in the proof
of part (i), so
$\partial(x^*)=-(\partial(x))^*$.
\hfill $\Box$ \\

\begin{lemma}
\label{lem:dr}
For $n \geq 1$
 the kernel of $\partial$ on $A_n$
is $A_n\cap A^*_n$.
\end{lemma}
\noindent {\it Proof:} 
 For $x \in A_n$ we show that
$\partial(x)=0$ if and only if $x \in A^*_n$.
First assume $\partial(x)=0$. Then 
$x^*=(-1)^{n-1}x$ by
(\ref{eq:par}), so $x\in A^*_n$.
To get the reverse implication, assume
$x \in A^*_n$ and note that $x^*\in A_n$.
Now each of $x, x^*$ is contained in $A_n$, so
$\partial(x) \in A_n$ in view of
(\ref{eq:par}). 
But $\partial(x) \in A_{n+1}$ by Lemma
\ref{lem:parset} and $A_n \cap A_{n+1}=0$
by Theorem
\ref{lem:ds} so $\partial(x)=0$.
\hfill $\Box$ \\

\noindent Our next goal is to show that
for $n\geq 1$ the image of $A_n$ under $\partial$
is $A_{n+1}\cap A^*_{n+1}$.
To this end it will be convenient to introduce
some subspaces
${}^+A_n$ and
${}^0A_n$ of $A_n$.

\begin{definition}
\label{def:tnpos}
\rm For $n\geq 1$
let 
${}^+A_n$
(resp. ${}^0A_n$)
denote the subspace of $A_n$ with a basis
consisting of the NR words that
have length $n$, and end with one of
$e_1, e_2, \ldots, e_d$
(resp. end with $e_0$).
\end{definition}

\begin{example}
\rm
Assume $d=2$. The basis for 
${}^+A_3$ from Definition
\ref{def:tnpos} is
\begin{eqnarray*}
&&e_1e^*_0e_1,
\qquad 
e_2e^*_0e_1,
\qquad e_0e^*_2e_1,
\qquad e_1e^*_2e_1,
\\
&&e_1e^*_0e_2,
\qquad 
e_2e^*_0e_2,
\qquad e_0e^*_1e_2,
\qquad e_2e^*_1e_2.
\end{eqnarray*}
\noindent The basis for 
${}^0A_3$ from Definition
\ref{def:tnpos} is
\begin{eqnarray*}
&&e_0e^*_1e_0,
\qquad 
e_2e^*_1e_0,
\qquad e_0e^*_2e_0,
\qquad e_1e^*_2e_0.
\end{eqnarray*}
\end{example}

\begin{lemma}
\label{lem:tposobv}
For $n\geq 1$,
\begin{itemize}
\item[\rm (i)] $A_n = {}^+A_n + {}^0A_n$ {\rm (direct sum)}.
\item[\rm (ii)] The dimension of  $\,{}^+A_n$ is $d^n$.
\item[\rm (iii)] The dimension of  $\,{}^0A_n$ is $d^{n-1}$.
\end{itemize}
\end{lemma} 
\noindent {\it Proof:} 
Routine using
Lemma
\ref{lem:sd} and
Definition
\ref{def:tnpos}.
\hfill $\Box$ \\

\begin{definition}
\label{def:sigmap}
\rm
For $n\geq 1$ we define an isomorphism
of vector spaces
$\sigma :
 {}^+A_n \to {}^0A_{n+1}$.
To do this we give the action of $\sigma$
on the basis for 
 ${}^+A_n$
from Definition
\ref{def:tnpos}.
Let $(r_1, r_2, \ldots, r_n)$
denote an NR word in $A$ such that $r_n\not=0$.
We  define the image of this word under
$\sigma$ to be
 $(r_1, r_2, \ldots, r_n,0)$. Note that
$\sigma$ sends the above basis for 
 ${}^+A_n$
to the basis
for 
 ${}^0A_{n+1}$ given in
Definition
\ref{def:tnpos}.
 Therefore $\sigma$ is an isomorphism of
vector spaces.
\end{definition} 

\begin{lemma}
\label{lem:sigalt}
For $n\geq 1$ and $x \in 
 {}^+A_n$,
\begin{eqnarray}
\sigma(x) = (-1)^n \partial(x)e_0.
\label{eq:sigalt}
\end{eqnarray}
\end{lemma}
\noindent {\it Proof:} 
Without loss we may assume
that $x$ is a vector in the basis for
 ${}^+A_n$ given in Definition
\ref{def:tnpos}.
Thus $x$ is an NR word $(r_1, r_2, \ldots, r_n)$
such that $r_n\not=0$.
Observe that $xe_0=0$ and
$x^*e_0=(r_1, r_2, \ldots, r_n,0)$.
By this and
(\ref{eq:par}) we find
$(-1)^n\partial(x)e_0$ is equal to 
$(r_1, r_2, \ldots, r_n,0)$, which is equal to
$\sigma(x)$ by Definition
\ref{def:sigmap}.
The result follows.
\hfill $\Box$ \\

\begin{lemma}
\label{lem:sigds}
For $n\geq 1$,
\begin{eqnarray}
A_n = {}^+A_n + A_n \cap A^*_n \qquad \hbox{\rm (direct sum)}.
\label{eq:sigds}
\end{eqnarray}
Moreover the dimension of $A_n \cap A^*_n$ is $d^{n-1}$.
\end{lemma}
\noindent {\it Proof:} 
We first show that the sum
${}^+A_n + A_n \cap A^*_n$
is direct. By Lemma
\ref{lem:sigalt} and since
 $\sigma :{}^+A_n\to {}^0A_{n+1}$ is
a bijection, the restriction of
$\partial$ to  
${}^+A_n$
is injective. Therefore the kernel
of $\partial $ on 
$A_n$ has 
  zero intersection with
${}^+A_n$.
This kernel is
$A_n \cap A^*_n$ by
Lemma
\ref{lem:dr}. Therefore
${}^+A_n$ has zero intersection 
with
$A_n \cap A^*_n$ so the sum
${}^+A_n + A_n \cap A^*_n$
is direct.
Let $k_n$ denote the dimension of
$A_n \cap A^*_n$. 
By our comments so far,
and given the dimensions of
$A_n$ and ${}^+A_n$ 
from Lemma
\ref{lem:dimtn}
and
Lemma
\ref{lem:tposobv}, respectively,
we obtain $k_n \leq d^{n-1}$, with equality if and only if
$A_n = {}^+A_n + A_n \cap A^*_n$.
To finish the proof it suffices to show
$k_n=d^{n-1}$. We do this by induction on $n$.
First assume $n=1$. We have $k_1 \leq 1$
by our above remarks, and $k_1\geq 1$
since $1 \in A_1 \cap A^*_1$ by
(\ref{eq:rel2}). Therefore $k_1=1$ as desired.
Next assume $n\geq 2$. Let $I_n$ denote the
image of $A_{n-1}$ under $\partial$.
By linear algebra the dimension of
$I_n$ is equal to the dimension of $A_{n-1}$
minus the dimension of the kernel of $\partial$
on $A_{n-1}$. The dimension of
$A_{n-1}$ is $d^{n-1}+d^{n-2}$.
The kernel of $\partial $ on
$A_{n-1}$ is $A_{n-1}\cap A^*_{n-1}$
so its dimension is $k_{n-1}$, which is
$d^{n-2}$ by induction.
Therefore the dimension of $I_n$ is
$d^{n-1}$. By 
Theorem
\ref{thm:interp2} and
(\ref{eq:par})
we have 
$I_n \subseteq 
A_n\cap A^*_n$. In this inclusion we consider
the dimensions and get
$d^{n-1}\leq k_n$. We showed earlier that
$k_n \leq d^{n-1}$ so 
$k_n =d^{n-1}$ as desired. The result follows.
\hfill $\Box$ \\

\begin{lemma}
\label{lem:image}
For $n \geq 1$
the image of $A_n$ under $\partial$ is
 $A_{n+1}\cap A^*_{n+1}$.
\end{lemma}
\noindent {\it Proof:} 
Denote this image by 
$I_{n+1}$, and
observe 
$I_{n+1} \subseteq 
 A_{n+1}\cap A^*_{n+1}$ by
Theorem
\ref{thm:interp2}.  To finish
the proof we show that
$I_{n+1}$ and
$ A_{n+1}\cap A^*_{n+1}$ have
the same dimension. By
Lemma
\ref{lem:sigds} the dimension of
$ A_{n+1}\cap A^*_{n+1}$ is $d^n$.
By Lemma
\ref{lem:dr} and 
(\ref{eq:sigds}) 
the dimension of 
$I_{n+1}$ is equal to
the dimension of
${}^+A_n$, which is $d^n$ by
Lemma
\ref{lem:tposobv}(ii).
The result follows. 
\hfill $\Box$ \\

\begin{lemma}
\label{lem:intone}
We have $A_1 \cap A^*_1 = \F 1$.
\end{lemma}
\noindent {\it Proof:} 
Observe
$\F 1 \subseteq A_1 \cap A^*_1$
by (\ref{eq:rel2}), and 
$A_1 \cap A^*_1$ has dimension 1 by
Lemma
\ref{lem:sigds}. 
\hfill $\Box$ \\

\begin{definition}
\rm
\label{def:incl}
Let $\iota :\F \to A$ denote $\F$-algebra
homomorphism that sends $a \mapsto a 1$ for $a \in \F$.
Note that 
 $\iota$ is an injection.
\end{definition}

\begin{theorem}
\label{thm:exact} The sequence
\[
  \begin{CD}
  \F 
 @>>\iota>
          A_1 
 @>>\partial>
          A_2
 @>>\partial>
          A_3
 @>>\partial>
      \cdots 
 \\
  \end{CD}
\]
is exact in the sense of
{\rm \cite[p.~435]{rotman}}.

\end{theorem}
\noindent {\it Proof:} 
This follows from
Lemma
\ref{lem:dr},
Lemma
\ref{lem:image},
 and
Lemma
\ref{lem:intone}.
\hfill $\Box$ \\

\noindent We emphasize a few points for later use.

\begin{lemma}
\label{lem:isoemph}
For $n\geq 1$ the restriction of
$\partial $ to
${}^+A_n$
is an isomorphism of vector spaces
${}^+A_n \to A_{n+1}\cap A^*_{n+1}$.
\end{lemma}
\noindent {\it Proof:} 
Combine Lemma
\ref{lem:dr},
line (\ref{eq:sigds}),
and Lemma
\ref{lem:image}.
\hfill $\Box$ \\

\begin{lemma}
\label{lem:pol}
For $n\geq 1$ and $x \in A_n$ the following
are equivalent:
\begin{itemize}
\item[\rm (i)] 
$x^*= (-1)^{n-1}x$;
\item[\rm (ii)] 
$x \in A_n \cap A^*_n$.
\end{itemize}
\end{lemma}
\noindent {\it Proof:} 
Combine 
(\ref{eq:par})
and 
Lemma
\ref{lem:dr}.
\hfill $\Box$ \\

\begin{lemma}
For $n\geq 1$ the map $\partial $ acts on
$A^*_n$ as follows. 
\begin{eqnarray*}
\partial(y) = -y-(-1)^n y^* \qquad \qquad (\forall y \in A^*_n).
\end{eqnarray*}
\end{lemma}
\noindent {\it Proof:} 
Write $x=y^*$, so that
$x \in A_n$ and
$y=x^*$.
Now compute $\partial(y)$ using
Lemma
\ref{lem:bound}(ii) and
 (\ref{eq:par}).
\hfill $\Box$ \\


\section{A subalgebra of $A$}

\noindent In this section we consider the sum
\begin{eqnarray}
\label{eq:tsum}
\sum_{n=0}^\infty 
(A_{n+1}\cap A^*_{n+1}).
\end{eqnarray}
We observe by Theorem
\ref{lem:ds}
and 
Lemma
\ref{lem:dr}
that
(\ref{eq:tsum})
is the kernel of
the map $\partial : A \to A$.
We will show that 
(\ref{eq:tsum})
is 
a subalgebra of $A$ that is
free of rank $d$.

\begin{lemma}
\label{lem:grading}
For nonnegative integers $n,m$ 
the following {\rm (i)--(iii)} hold.
\begin{itemize}
\item[\rm (i)] 
$
(A_{n+1}\cap A^*_{n+1})A_{m+1} \subseteq A_{n+m+1}
$;
\item[\rm (ii)]
$
(A_{n+1}\cap A^*_{n+1})A^*_{m+1} \subseteq A^*_{n+m+1}
$;
\item[\rm (iii)]
$
(A_{n+1}\cap A^*_{n+1})
(A_{m+1}\cap A^*_{m+1})
\subseteq
A_{n+m+1}\cap A^*_{n+m+1}
$.
\end{itemize}
\end{lemma}
\noindent {\it Proof:} 
(i)
For $x \in 
A_{n+1}\cap A^*_{n+1}$
and
 $y \in 
A_{m+1}$
we show that
$xy \in 
A_{n+m+1}$.
First assume $m$ is even.
Using 
 $x \in 
A_{n+1}$
and
 $y \in 
A_{m+1}$
and the inclusion on the left in
(\ref{eq:tntmodd}), we obtain 
$xy \in 
A_{n+m+1}
$.
Next assume $m$ is odd.
Using 
 $x \in 
A^*_{n+1}$
and
 $y \in 
A_{m+1}$
and the inclusion on the right in
(\ref{eq:tntmeven}), we obtain
$xy \in 
A_{n+m+1}
$.
\\
\noindent (ii)
For $x \in 
A_{n+1}\cap A^*_{n+1}$
and
 $y \in 
A^*_{m+1}$
we show that
$xy \in 
A^*_{n+m+1}$.
Observe that
$x^* \in 
A_{n+1}\cap A^*_{n+1}$
and 
$y^* \in 
A_{m+1}$
so 
$x^*y^* \in A_{n+m+1}$
by (i) above.
Applying $*$ we find
$xy \in A^*_{n+m+1}$.
\\
\noindent (iii)
Combine (i) and (ii) above.
\hfill $\Box$ \\

\begin{corollary}
\label{cor:subalg}
The sum
{\rm (\ref{eq:tsum})}
is a subalgebra of $A$.
\end{corollary}
\noindent {\it Proof:} 
The sum contains the identity $1$ of $A$ by
Lemma
\ref{lem:intone}. The sum
is closed under multiplication
by Lemma
\ref{lem:grading}(iii).
\hfill $\Box$ \\

\noindent We will return to the subalgebra
(\ref{eq:tsum})
after a few comments.

\begin{lemma}
\label{lem:altform}
For each basis vector $(r_1, r_2, \ldots, r_n)$ 
from Theorem
\ref{thm:tbasis}(i), the element
\begin{eqnarray*}
(r_1, r_2, \ldots, r_n)
\quad + \quad 
(-1)^n(r_1, r_2, \ldots, r_n)^*
\end{eqnarray*}
is equal to
\begin{eqnarray}
(e_{r_1}-e^*_{r_1})
(e_{r_2}-e^*_{r_2})
\cdots
(e_{r_n}-e^*_{r_n})(-1)^{\lfloor n/2 \rfloor}.
\label{eq:redexpand}
\end{eqnarray}
The expression
${\lfloor x \rfloor}$ denotes the greatest integer less than
or equal to $x$.
\end{lemma}
\noindent {\it Proof:} 
Expand 
(\ref{eq:redexpand}) into a sum of $2^n$ terms.
Simplify these terms using
(\ref{eq:rel1}) and 
$r_{i-1}\not=r_i$ for $2 \leq i \leq n$.
\hfill $\Box$ \\

\medskip
\noindent Let $F$ denote the 
$\F$-algebra defined by
generators 
 $\lbrace s_i\rbrace_{i=1}^d$ and no relations.
Thus $F$ is the free $\F$-algebra of rank $d$.
We call $\lbrace s_i\rbrace_{i=1}^d$ the {\it 
standard generators} for $F$.
We recall a few facts about $F$. For an integer $n\geq 0$,
by a {\it word in $F$ of length $n$} we mean a product
$y_1y_2\cdots y_n$
such that
$\lbrace y_i\rbrace_{i=1}^n$
are 
standard generators for $F$.
We interpret the word of length $0$ to be the identity of $F$.
The $\F$-vector space $F$ has a basis consisting
of its words
\cite[p.~723]{rotman}.
For $n\geq 0$ let $F_n$ denote the subspace of $F$ spanned by
the words of length $n$. 
Note that $F_n$ has dimenion $d^n$.
We have a direct sum
$F=\sum_{n=0}^\infty F_n$, and $F_rF_s =F_{r+s}$
for $r,s\geq 0$. 
We call $F_n$
the {\it $n$th homogeneous component} of $F$.

\begin{theorem}
\label{thm:freeinj}
With the above notation, 
consider the $\F$-algebra homomorphism $F\to A$
that sends $s_i \mapsto e_i-e^*_i$ for $1 \leq i \leq d$.
This map is an injection
and its image is $\sum_{n=0}^\infty (A_{n+1}\cap A^*_{n+1})$.
Moreover for $n\geq 0$
the image of  $F_n$
is  $A_{n+1}\cap A^*_{n+1}$.
\end{theorem}
\noindent {\it Proof:} 
Let $\varepsilon : F\to A$ denote the homomorphism in question.
We claim that for $n\geq 0$
the restriction of $\varepsilon $ to
$F_n$ is a bijection
$F_n \to A_{n+1}\cap A^*_{n+1}$.
To establish the claim we split the argument
into three cases: $n=0$, $n=1$, 
and $n\geq 2$.
The claim holds for $n=0$ by
Lemma
\ref{lem:intone} and since
$F_0=\F 1$.
To see that the claim holds for $n=1$,
note
 that $F_1$ has a basis ${\lbrace s_i\rbrace}_{i=1}^d$.
By 
Definition
\ref{def:tnpos}
the elements $\lbrace e_i\rbrace_{i=1}^d$
form a basis for 
${}^+A_1$,
so $\lbrace \partial(e_i)\rbrace_{i=1}^d$
is a basis for $A_2 \cap A^*_2$ in 
view of
Lemma
\ref{lem:isoemph}.
By Lemma
\ref{lem:map}
 we have
$\partial(e_i)=e_i-e^*_i$
for $1 \leq i \leq d$.
Therefore 
$\lbrace e_i-e^*_i\rbrace_{i=1}^d$
is a basis for $A_2\cap A^*_2$, 
and the claim follows for $n=1$.
We now show that the claim holds for $n\geq 2$.
Using
$F_n=(F_1)^n$,
$\varepsilon(F_1)=A_2\cap A^*_2$,
 and 
Lemma
\ref{lem:grading}(iii) we obtain
$\varepsilon(F_n) \subseteq A_{n+1}\cap A^*_{n+1}$.
To see the reverse inclusion,
first note by Lemma
\ref{lem:image}
that any element in
$A_{n+1}\cap A^*_{n+1}$ can be written
as $\partial(x)$ for some $x \in A_n$.
We show  
$\partial(x) \in \varepsilon (F_n)$.
Without loss we may assume
that
$x$ is a vector
$(r_1, r_2, \ldots,r_n)$ in the
basis for $A_n$
from
Lemma
\ref{lem:sd}. 
Combining
Lemma
\ref{lem:map} and
Lemma
\ref{lem:altform}
we find that 
$\partial(x)$ is equal to
(\ref{eq:redexpand}). In particular
$\partial(x)=
\varepsilon(z_1)
\varepsilon(z_2)
\cdots
\varepsilon(z_n)
$
where $z_i \in F_1$
for $1 \leq i \leq n$.
Observe
$\partial(x)=\varepsilon(z_1z_2\cdots z_n)$
and $z_1z_2\cdots z_n \in F_n$
so 
 $\partial(x)\in \varepsilon(F_n)$.
Therefore
$A_{n+1}\cap A^*_{n+1}
\subseteq 
\varepsilon(F_n)$.
So far we have $\varepsilon(F_n)=A_{n+1}\cap A^*_{n+1}$.
To show that the map
$F_n \to A_{n+1}\cap A^*_{n+1}$, 
$x \mapsto \varepsilon(x)$
is a bijection, it suffices to 
show that
$F_n$ and $A_{n+1}\cap A^*_{n+1}$ have
the same dimension. We mentioned
below Lemma 
\ref{lem:altform} that
$F_n$ has dimension $d^n$.
By the last line of
Lemma
\ref{lem:sigds} we find
 $A_{n+1}\cap A^*_{n+1}$ 
also has dimension $d^n$.
By these 
comments
the map $F_n \to A_{n+1}\cap A^*_{n+1}$, 
$x \mapsto \varepsilon(x)$
is a bijection. The claim is now proved for $n\geq 2$.
We have established the claim,
and the result follows in view of
the directness of the
 sum
 $\sum_{n=0}^\infty A_{n+1} \cap A^*_{n+1}$.
\hfill $\Box$ \\
 
\noindent For notational convenience let us
identify the free algebra $F$ from above
Theorem 
\ref{thm:freeinj}
with
the subalgebra 
(\ref{eq:tsum}) of $A$,
via the injection
from 
Theorem
\ref{thm:freeinj}.
Our next goal is to show that each of the
$\F$-linear maps
\begin{eqnarray*}
&&F \otimes A_1 \;\to\; A 
\qquad \qquad \qquad 
F \otimes A^*_1 \; \;\to \;\; A
\\
&&u \otimes v \; \; \mapsto \; \; uv
\qquad \qquad \qquad 
u \otimes v \;  \;\mapsto \;  \; uv
\end{eqnarray*}
is an isomorphism of $\F$-vector spaces.
We need a lemma.

\begin{lemma}
\label{lem:part1}
For positive integers $n,m$ the 
$\F$-linear map
\begin{eqnarray*}
(A_n \cap A^*_n) \otimes A_m &\to& A_{n+m-1}
\\
u \otimes v \;  &\mapsto & \; uv
\end{eqnarray*}
is an isomorphism of $\F$-vector spaces.
\end{lemma}
\noindent {\it Proof:} 
Let
$\theta$
denote the map in question.
To show that 
$\theta$
is bijective, we show
that the dimension of
$(A_n \cap A^*_n) \otimes A_m$
is equal to the dimension of $A_{n+m-1}$,
and that
$\theta$
is surjective.
The dimension of 
$A_n \cap A^*_n$ is
$d^{n-1}$
by Lemma
\ref{lem:sigds}, and the dimension of
$A_m$ is
$(d+1)d^{m-1}$ by Lemma
\ref{lem:dimtn}, so the dimension of
$(A_n \cap A^*_n) \otimes A_m$
is $(d+1)d^{n+m-2}$.
The dimension of
 $A_{n+m-1}$
is $(d+1)d^{n+m-2}$ by
Lemma \ref{lem:dimtn}. Therefore the dimensions of
$(A_n \cap A^*_n) \otimes A_m$
and $A_{n+m-1}$ are the same.
Next we show that 
$\theta$
is surjective.
First assume $n=1$. Then 
$\theta$
is surjective
since $1 \in 
A_1 \cap A^*_1$ by
Lemma
\ref{lem:intone}.
Next assume $n\geq 2$. By
Lemma
\ref{lem:sd} the space
$A_{n+m-1}$ has a basis consisting
of the NR words in $A$ that have 
length $n+m-1$ and end with a nonstarred
element. We show that each of these basis
elements is in the image of $\theta $.
Consider an NR word $w=(r_1, r_2, \ldots, r_{n+m-1})$.
Define $u=(-1)^{(n-1)m}\partial(r_1, r_2, \ldots, r_{n-1})$
and observe that $u \in A_n\cap A^*_n$ by
Lemma
\ref{lem:image}.
Define $v=(r_n, r_{n+1}, \ldots, r_{n+m-1})$
and observe $v \in A_m$.
One verifies $w=uv$ 
by first using 
(\ref{eq:par}),
and then 
Theorem
\ref{thm:prod}(i),
Theorem
\ref{thm:prod2}(iii)
if $m$ is odd and
Theorem
\ref{thm:prod}(iii),
Theorem
\ref{thm:prod2}(i)
if $m$ is even.
 Therefore
$w$ is the image of $u\otimes v$ under
$\theta$.
We have shown 
that
$\theta$
is surjective, and 
the result follows.
\hfill $\Box$ \\

%

\begin{theorem}
Each of the $\F$-linear maps
\begin{eqnarray*}
&&F \otimes A_1 \;\to\; A 
\qquad \qquad \qquad 
F \otimes A^*_1 \; \;\to \;\; A
\\
&&u \otimes v \; \; \mapsto \; \; uv
\qquad \qquad \qquad 
u \otimes v \;  \;\mapsto \;  \; uv
\end{eqnarray*}
is an isomorphism of $\F$-vector spaces.
\end{theorem}
\noindent {\it Proof:} 
Let $\psi$ (resp. $\xi$) denote the map on
the left (resp. right).
We first show that $\psi$ is an isomorphism
of $\F$-vector spaces.
By construction
the sum
$F=\sum_{n=1}^\infty A_n \cap A^*_n$
is direct. Therefore
the sum
\begin{eqnarray*}
F \otimes A_1 = \sum_{n=1}^\infty (A_n \cap A^*_n)\otimes A_1
\end{eqnarray*}
is direct. By Theorem
\ref{lem:ds} the sum
$A = \sum_{n=1}^\infty A_n$
is direct. For $n \geq 1$ we apply Lemma
\ref{lem:part1}
with $m=1$ and find that
the map
\begin{eqnarray*}
&&(A_n \cap A^*_n) \otimes A_1 \;\;\to \;\; A_{n}
\\
&& \qquad \qquad u \otimes v \; \;\mapsto \; \; uv
\end{eqnarray*}
is an isomorphism of $\F$-vector spaces.
It follows that $\psi$ is
an isomorphism of $\F$-vector spaces.
The map $\xi$ is an isomorphism of
 $\F$-vector spaces since it is  
the
composition 

\[
  \begin{CD}
  F \otimes A^*_1 
 @>>* \otimes *>
  F \otimes A_1 
 @>>\psi>
         A 
 @>>*>
          A
  \end{CD}
\]
and each composition factor is an isomorphism
of 
$\F$-vector spaces.
\hfill $\Box$ \\

\section{The subspaces $A_n$ revisited}

\noindent In this section we
 present a more detailed version of
Theorem
\ref{thm:ptntm}. Let $n,m$ denote positive integers.
For $m$ odd we consider the $\F$-linear maps
\begin{eqnarray*}
&&A_n \otimes A_m \;\;\to \;\; A_{n+m-1}
\qquad \qquad \qquad 
A_n \otimes A_m \;\;\to\;\; A_{n+m}+ A_{n+m-1}
\\
&& \;\quad u \otimes v \;\;  \mapsto \; \; uv
\qquad \qquad \qquad 
\qquad \;\;\quad  
u \otimes v \; \; \mapsto \; \; u^*v
\end{eqnarray*}
and for $m$ even we consider the $\F$-linear maps
\begin{eqnarray*}
&&A_n \otimes A_m \;\;\to\;\; A_{n+m}+ A_{n+m-1}
\qquad \qquad \qquad 
A_n \otimes A_m \;\;\to\;\;  A_{n+m-1}
\\
&&\quad \;u \otimes v \;  \;\mapsto \; \; uv
\qquad \qquad \qquad 
\qquad  \qquad \qquad 
\;u \otimes v \;  \;\mapsto \; \; u^*v.
\end{eqnarray*}

\begin{definition}
\label{def:eqneq}
\rm
Let $n,m$ denote positive integers.
\begin{itemize}
\item[\rm (i)]
Let
 ${}^{\not=}
(A_n\otimes A_m)$
denote the subspace of 
  $A_n \otimes A_m$ that has a basis
consisting
of the elements $u \otimes v$, where 
$u=(r_1, r_2, \ldots, r_n)$ is an NR word in $A_n$
and $v=(r'_1, r'_2, \ldots, r'_m)$ is an NR word in $A_m$ such that
$r_n \not=r'_1$.
\item[\rm (ii)]
Let
 ${}^=
(A_n\otimes A_m)$
denote the subspace of 
  $A_n \otimes A_m$ that has a basis
consisting
of the elements $u \otimes v$, where 
$u=(r_1, r_2, \ldots, r_n)$ is an NR word in $A_n$
and $v=(r'_1, r'_2, \ldots, r'_m)$ is an NR word in $A_m$ such that
$r_n =r'_1$.
\end{itemize}
\end{definition}

\noindent The following result is immediate
from Definition
\ref{def:eqneq}.

\begin{lemma}
\label{lem:basic}
With reference to Definition
\ref{def:eqneq},
\begin{eqnarray}
A_n \otimes A_m &=& 
 {}^{\not=}
(A_n\otimes A_m)
\quad +  \quad
 {}^=
(A_n\otimes A_m) \qquad \qquad \mbox{\rm (direct sum)}.
\end{eqnarray}
\end{lemma}

\begin{theorem}
\label{thm:kerok}
For positive integers $n,m$ the following
{\rm (i), (ii)} hold.
\begin{itemize}
\item[\rm (i)]
Assume $m$ is odd.
Then the $\F$-linear map
\begin{eqnarray*}
A_n \otimes A_m &\to& A_{n+m-1}
\\
u \otimes v \;  &\mapsto & \; uv
\end{eqnarray*}
is surjective with kernel
 $\;{}^{\not=}
(A_n\otimes A_m)$.
\item[\rm (ii)]
Assume $m$ is even.
Then the $\F$-linear map
\begin{eqnarray*}
A_n \otimes A_m &\to& A_{n+m-1}
\\
u \otimes v \;  &\mapsto & \; u^*v
\end{eqnarray*}
is surjective with kernel
 $\;{}^{\not=}
(A_n\otimes A_m)$.
\end{itemize}
\end{theorem}
\noindent {\it Proof:} 
(i) A basis for
 $\;{}^{\not=}
(A_n\otimes A_m)$ is given in
Definition
\ref{def:eqneq}(i).
By Theorem
\ref{thm:prod}(i)
the map sends each element in this
basis to zero.
A basis for
 $\;{}^{=}
(A_n\otimes A_m)$ is given in
Definition
\ref{def:eqneq}(ii).
By Theorem
\ref{thm:prod}(ii)
the map sends 
this basis 
to the basis for $A_{n+m-1}$ given in
Lemma
\ref{lem:sd}.
The result follows from these comments
and Lemma
\ref{lem:basic}.
\\
\noindent (ii) Similar to the proof of 
(i) above.
\hfill $\Box$ \\

\begin{lemma}
\label{prop:concl}
For positive integers $n,m$ we have
\begin{eqnarray}
A_n \otimes A_m &=& 
 {}^{\not=}
(A_n\otimes A_m)
\quad + \quad 
(A_n\cap A^*_n)\otimes A_m \qquad \qquad \mbox{\rm (direct sum).}
\end{eqnarray}
\end{lemma}
\noindent {\it Proof:} 
For $m$ odd the result follows from
Lemma
\ref{lem:part1} and 
Theorem
\ref{thm:kerok}(i). For $m$ even
the result follows from
Lemma
\ref{lem:pol},
Lemma
\ref{lem:part1},
and 
Theorem
\ref{thm:kerok}(ii).
\hfill $\Box$ \\

\noindent The following result will be helpful.

\begin{proposition}
\label{prop:partialmap}
For positive integers $n,m$
the $\F$-linear map
\begin{eqnarray*}
A_n \otimes A_m &\to& A_{n+m}
\\
u \otimes v \;  &\mapsto & \; \partial(u)v
\end{eqnarray*}
is surjective 
with  kernel
$(A_n \cap A^*_n)\otimes A_m$.
\end{proposition}
\noindent {\it Proof:} 
By Lemma
\ref{lem:dr} and Lemma
\ref{lem:image}, the map
$A_n \to A_{n+1}\cap A^*_{n+1}$,
$u \mapsto \partial(u)$
is surjective with kernel
$A_n \cap A^*_n$.
Therefore the map
$A_n \otimes A_m \to (A_{n+1}\cap A^*_{n+1})\otimes A_m$,
$u\otimes v \mapsto \partial(u)\otimes v$
is surjective with kernel
$(A_n \cap A^*_n)\otimes A_m$.
By Lemma
\ref{lem:part1}
the map
$(A_{n+1}\cap A^*_{n+1})\otimes A_m \to A_{n+m}$,
$u \otimes v\mapsto uv$
is a bijection. Composing the two previous maps,
we find that the map
$A_n \otimes A_m \to A_{n+m}$, $u \otimes v \mapsto \partial(u)v$
is surjective with kernel
$(A_n\cap A^*_n)\otimes A_m$.
\hfill $\Box$ \\

\begin{theorem}
For positive integers $n,m$ the following {\rm (i), (ii)} hold.
\begin{itemize}
\item[\rm (i)] Assume $m$ is even. Then 
the $\F$-linear map
\begin{eqnarray*}
A_n \otimes A_m &\to& A_{n+m}+ A_{n+m-1}
\\
u \otimes v \;  &\mapsto & \; uv
\end{eqnarray*}
is an isomorphism of $\F$-vector spaces.
Under this map the preimage of
$A_{n+m}$ is
 $\;{}^{\not=}
(A_n\otimes A_m)$ and the preimage of
$A_{n+m-1}$ is
$(A_n \cap A^*_n)\otimes A_m$.
\item[\rm (ii)]
Assume $m$ is odd. Then 
the $\F$-linear map
\begin{eqnarray*}
A_n \otimes A_m &\to& A_{n+m}+ A_{n+m-1}
\\
u \otimes v \;  &\mapsto & \; u^*v
\end{eqnarray*}
is an isomorphism of $\F$-vector spaces.
Under this map the preimage of
$A_{n+m}$ is
 $\;{}^{\not=}
(A_n\otimes A_m)$ and the preimage of
$A_{n+m-1}$ is
$(A_n \cap A^*_n)\otimes A_m$.
\end{itemize}
\end{theorem}
\noindent {\it Proof:} 
For $u \in A_n$ and $v \in A_m$
we use $\partial(u)=u+(-1)^n u^*$ 
to obtain
\begin{eqnarray}
\label{eq:breakapart}
\partial(u)v= uv+(-1)^n u^*v.
\end{eqnarray}
\noindent (i) Denote the map by $\eta $.
The restriction of $\eta$ to
 $\;{}^{\not=}
(A_n\otimes A_m)$
gives a bijection
 $\;{}^{\not=}
(A_n\otimes A_m) \to A_{n+m}$
by Theorem
\ref{thm:kerok}(ii), 
Lemma
\ref{prop:concl},
Proposition
\ref{prop:partialmap}, and
(\ref{eq:breakapart}).
The restriction of
$\eta $ to
$(A_n \cap A^*_n)\otimes A_m$
gives a bijection
$(A_n \cap A^*_n)\otimes A_m \to A_{n+m-1}$,
by 
Lemma \ref{lem:part1}.
The result follows.
\\
\noindent (ii) Denote the map by $\zeta $.
The restriction of $\zeta$ to
 $\;{}^{\not=}
(A_n\otimes A_m)$
gives a bijection
 $\;{}^{\not=}
(A_n\otimes A_m) \to A_{n+m}$
by Theorem
\ref{thm:kerok}(i), 
Lemma
\ref{prop:concl},
Proposition
\ref{prop:partialmap}, and
(\ref{eq:breakapart}).
The restriction of
$\zeta$ to
$(A_n \cap A^*_n)\otimes A_m$
gives a bijection
$(A_n \cap A^*_n)\otimes A_m \to A_{n+m-1}$,
by 
Lemma \ref{lem:pol}
and
Lemma \ref{lem:part1}.
The result follows.
\hfill $\Box$ \\

\section{The subspaces $A_{\leq n}$}

\noindent In this last section we investigate the
following subspaces of $A$.

\begin{definition}
\label{def:leqn}
\rm
For all integers $n\geq 1$ we define
\begin{eqnarray}
 A_{\leq n} = A_1 + A_2 + \cdots + A_n.
\label{eq:defleqn}
\end{eqnarray}
\end{definition}

\noindent The following lemma is immediate from the construction.

\begin{lemma}
\label{lem:obv}
For  $n\geq 1$ we display a basis for each relative of
$A_{\leq n}$.
\begin{center}
\begin{tabular}{c|c}
space  &  basis \\
\hline
$A_{\leq n}$ &
the NR words in $A$ that have length at most
 $n$ and end with a nonstarred element
\\
$A^*_{\leq n}$ &
the NR words in $A$ that have length at most $n$ and end with a starred element
\\
$A^\dagger_{\leq n}$ &
the NR words in $A$ that have length at most $n$ and begin with a nonstarred element
\\
$A^{*\dagger}_{\leq n}$ &
the NR words in $A$ 
that have length at most $n$ and begin with a starred element
\end{tabular}
\end{center}
\end{lemma}

\begin{lemma}
\label{lem:dimleqn}
For $n\geq 1$ the relatives of 
$A_{\leq n}$ all have dimension
$(d+1)(1+d+d^2+\cdots+d^{n-1})$.
\end{lemma}
\noindent {\it Proof:} 
By Theorem
\ref{lem:ds}
and Definition
\ref{def:leqn},
the dimension of
 $A_{\leq n}$ is equal to the sum of the dimensions of
$A_1, A_2, \ldots, A_n$.
The result follows from this and
Lemma
 \ref{lem:dimtn}.
\hfill $\Box$ \\

\noindent In Lemma
\ref{lem:obv} we gave a basis for each relative
of $A_{\leq n}$. 
In a moment we will display another basis.
In order to motivate this new basis we first give
a spanning set.

\begin{lemma}
\label{thm:leqn}
 For $n\geq 1$ we display a spanning
set for each relative of 
 $A_{\leq n}$.
\begin{center}
\begin{tabular}{c|c}
space  &  spanning set \\
\hline
$A_{\leq n}$ &
the words in $A$ that have length
 $n$ and end with a nonstarred element
\\
$A^*_{\leq n}$ &
the words in $A$ that have length  $n$ and end with a starred element
\\
$A^\dagger_{\leq n}$ &
the words in $A$ that have length  $n$ and begin with a nonstarred element
\\
$A^{*\dagger}_{\leq n}$ &
the words in $A$ that have length $n$ and begin with a starred element
\end{tabular}
\end{center}
\end{lemma}
\noindent {\it Proof:} 
Concerning the first row of the table, let $S_n$ denote
the subspace of $A$ spanned by the words in $A$
that have length $n$ and end with a nonstarred element.
We show $S_n=A_{\leq n}$.
By construction $S_n=\cdots A_1 A^*_1 A_1$ ($n$ factors).
By Theorem
\ref{thm:ptntm} 
we have
$A_1A_j \subseteq A_j+A_{j+1}$
and
$A^*_1A_j \subseteq A_j+A_{j+1}$
for $1 \leq j \leq n-1$. By this
and 
 induction on $n$
we find $S_n \subseteq A_{\leq n}$.
To get the reverse inclusion,
note that for $1 \leq j \leq n$ we have $A_j \subseteq S_j$,
and also $S_j \subseteq S_n$ 
since $ 1 \in A_1$ and $1 \in A^*_1$ by
Lemma
\ref{lem:intone}.
We have verified the first row of the table.
The remaining rows are similarly verified.
\hfill $\Box$ \\

\noindent 
For each spanning set in 
Lemma
\ref{thm:leqn}, the 
set is not a basis for $n\geq 3$,
since the set has cardinality $(d+1)^n$
and this number differs from the dimension
 given in
Lemma
\ref{lem:dimleqn}. Our next goal is
to obtain
a subset of the spanning set that is a basis.

\begin{definition}
\label{def:r/nr}
\rm
A word $g_1g_2\cdots g_n$ in $A$ is called
{\it repeating/nonrepeating} (or {\it R/NR})
whenever for $2 \leq j \leq n$,
if $g_{j-1}$, $g_j$ have the same index then
$g_1, g_2, \ldots, g_j$ all have the same index.
\end{definition}

\begin{example}
\rm
For $d=2$ we display the R/NR words in $A$ that have length 3
and end with $e_0$.
\begin{eqnarray*}
&&e_0e^*_0e_0, \qquad 
e_1e^*_1e_0, \qquad 
e_2e^*_2e_0, \qquad 
\\
&&
e_0e^*_1e_0, \qquad
e_2e^*_1e_0, \qquad
e_0e^*_2e_0, \qquad
e_1e^*_2e_0.
\end{eqnarray*}
\end{example}

\begin{definition}
\rm
A word in $A$ is called 
{\it nonrepeating/repeating} (or {\it NR/R}) whenever
its image under $\dagger$ is R/NR.
\end{definition} 

\begin{example}
\rm
For $d=2$ we display the NR/R words in $A$ that have length 3
and start with $e_0$.
\begin{eqnarray*}
&&e_0e^*_0e_0, \qquad 
e_0e^*_1e_1, \qquad 
e_0e^*_2e_2, \qquad 
\\
&&
e_0e^*_1e_0, \qquad
e_0e^*_1e_2, \qquad
e_0e^*_2e_0, \qquad
e_0e^*_2e_1.
\end{eqnarray*}
\end{example}

\begin{theorem} 
For $n\geq 1$ we display a basis for each relative of
$A_{\leq n}$.
\begin{center}
\begin{tabular}{c|c}
space  &  basis \\
\hline
$A_{\leq n}$ &
the R/NR words in $A$ that have length $n$ and end with a nonstarred element
\\
$A^*_{\leq n}$ &
the R/NR words in $A$ that have length $n$ and end with a starred element
\\
$A^\dagger_{\leq n}$ &
the NR/R words in $A$ that have length $n$ and begin with a nonstarred element
\\
$A^{*\dagger}_{\leq n}$ &
the NR/R words in $A$ that have length $n$ and begin with a starred element
\end{tabular}
\end{center}
\end{theorem}
\noindent {\it Proof:} 
Concerning the first row of the table,
let $(R/NR)_n$ denote the set of
R/NR words in $A$ that have length $n$ and end
with a nonstarred element. We show
$(R/NR)_n$ is a basis for $A_{\leq n}$.
Let $(NR)_n$ denote the basis
for $A_n$ given in Lemma
\ref{lem:sd}.
Let 
 $(NR)_{\leq n}=
\cup_{j=1}^n
 (NR)_j$
 denote the basis
for 
 $A_{\leq n}$ given in 
Lemma
\ref{lem:obv}.
We now define a linear transformation
$f:
 A_{\leq n} \to
 A_{\leq n}$.
To this end we give the action
of $f$ on 
 $(NR)_j$ for $1 \leq j\leq n$.
For a word $(r_1, r_2, \ldots, r_j)$ in
 $(NR)_j$
we define its  image under $f$ to be
$(r_1, r_1, \ldots,r_1, r_1, r_2, \ldots, r_j)$
($n$ coordinates). This image is contained in
$A_{\leq n}$ by
Lemma
\ref{thm:leqn}.
By the construction $f$ sends the
basis
$(NR)_{\leq n}$ to the set
$(R/NR)_n$.
To show that 
$(R/NR)_n$
 is a basis for $A_{\leq n}$ it suffices
to show that $f$ is a bijection.
Using the data in Theorem
\ref{thm:prod} and
Theorem
\ref{thm:prod2},
one finds
$(f-I)A_j \subseteq A_{j+1}+ \cdots+ A_n$ 
for $1 \leq j \leq n$, where $I:A\to A$ is the identity map.
 Therefore, with
respect to an appropriate ordering of the
basis 
$(NR)_{\leq n}$, the matrix
which represents $f$ is lower triangular, with
all diagonal entries $1$. This matrix is nonsingular
so $f$ is invertible and hence a bijection. Therefore
$(R/NR)_n$ is a basis
for 
 $A_{\leq n}$.
This yields the first row of the table. The remaining
rows are similarly obtained.
\hfill $\Box$ \\

\bigskip


\noindent Tatsuro Ito \hfil\break
\noindent Division of Mathematical and Physical Sciences \hfil\break
\noindent Graduate School of Natural Science and Technology\hfil\break
\noindent Kanazawa University \hfil\break
\noindent Kakuma-machi,  Kanazawa 920-1192, Japan \hfil\break
\noindent email:  {\tt tatsuro@kenroku.kanazawa-u.ac.jp}

\bigskip

\noindent Paul Terwilliger \hfil\break
\noindent Department of Mathematics \hfil\break
\noindent University of Wisconsin \hfil\break
\noindent 480 Lincoln Drive \hfil\break
\noindent Madison, WI 53706-1388 USA \hfil\break
\noindent email: {\tt terwilli@math.wisc.edu }\hfil\break

\end{document}